\documentclass[review, 11pt]{elsarticle}

\biboptions{authoryear}

\usepackage{hyperref}
\usepackage{rotating}

\hypersetup{
    colorlinks=true,
    linkcolor=blue,
    filecolor=magenta,      
    urlcolor=cyan,
    pdftitle={Overleaf Example},
    pdfpagemode=FullScreen,
    }
    
\usepackage[english]{babel}
\usepackage{amsmath,amsthm,amsfonts,amssymb,latexsym,amstext}

\usepackage[dvipsnames]{xcolor} 
\usepackage{color,soul}
\usepackage{mathrsfs} 
\usepackage{tabularx} 
\usepackage{arydshln}
\usepackage{natbib}
\usepackage[numbib]{tocbibind}
\usepackage{float} 
\usepackage{url}
\usepackage{multirow}
\newcolumntype{C}{>{\centering\arraybackslash}X}
\usepackage[skip=.5\baselineskip plus2pt, indent=20pt]{parskip} 
\usepackage{makecell}

\usepackage{ stmaryrd }
\usepackage{booktabs} 

\newtheoremstyle{thmstyle} 
{\topsep}                    
{\topsep}                    
{\itshape}                   
{}                           
{\bfseries}                   
{.}                          
{.5em}                       
{}  
\theoremstyle{thmstyle}
\newtheorem{theorem}{Theorem}[section]
\newtheorem{proposition}{Proposition}[section]

\theoremstyle{definition}
\newtheorem{definition}{Definition}[section]

\newtheorem{example}{Example}[section]

\usepackage{multirow}

\newcommand{\addQEDstyle}[2]{\AtBeginEnvironment{#1}{\pushQED{\qed}\renewcommand{\qedsymbol}{#2}}\AtEndEnvironment{#1}{\popQED}}
\addQEDstyle{example}{$\triangle$}

\newcommand\restr[2]{{
  \left.\kern-\nulldelimiterspace 
  #1 
  \vphantom{\big|} 
  \right|_{#2} 
  }}
  
\usepackage{mathtools} 

\usepackage[margin=1.5cm]{geometry}

 \usepackage{setspace} 
\makeatletter
\newcommand\thefontsize{The current font size is: \f@size pt}
\makeatother

\begin{document}

\begin{frontmatter}
\title{Highway toll \textcolor{black}{allocation} problem revisited: new methods and characterizations}

\author[mymainaddress]{Paula Soto-Rodr{\'i}guez}
\ead{paula.soto.rodriguez@rai.usc.es}

\author[mysecondaryaddress]{Balbina Casas-Méndez}
\ead{balbina.casas.mendez@usc.es}

\author[mythirdaddress]{Alejandro Saavedra-Nieves}
\ead{alejandro.saavedra.nieves@usc.es}

\address[mymainaddress]{\textcolor{black}{Corresponding author.} MODESTYA Research Group, Department of Statistics, Mathematical Analysis and Optimization, Faculty of Mathematics, Universidade de Santiago de Compostela, Campus Vida, 15782 Santiago de Compostela, Spain. ORCID: \href{https://orcid.org/0009-0005-5199-8223}{0009-0005-5199-8223}.}

\address[mysecondaryaddress]{CITMAga, MODESTYA Research Group, 
Department of Statistics, Mathematical Analysis and Optimization, Faculty of Mathematics, Universidade de Santiago de Compostela, Campus Vida, 15782 Santiago de Compostela, Spain. ORCID: \href{https://orcid.org/0000-0002-2826-218X}{0000-0002-2826-218X}.} 

\address[mythirdaddress]{CITMAga, MODESTYA Research Group, Department of Statistics, Mathematical Analysis and Optimization, Faculty of Mathematics, Universidade de Santiago de Compostela, Campus Vida, 15782 Santiago de Compostela, Spain. ORCID: \href{https://orcid.org/0000-0003-1251-6525}{0000-0003-1251-6525}.}

\begin{abstract} 
\textcolor{black}{This paper considers the highway toll allocation problem (Wu, van den Brink, and Est\'evez-Fern\'andez in Transport Res B-Meth 180:10288, 2024). The aim is to allocate the tolls collected from the users of a highway across the various road sections. To this end,} the authors propose, among others, the Segments Equal Sharing method, which is characterized and reinterpreted as a specific solution of a cooperative game associated with the problem. This paper presents two new allocation rules: the
Segments Proportional Sharing method 
and the 
Segments Compensated Sharing method. 
\textcolor{black}{We axiomatically characterize these new methods and compare their properties to those of the Segments Equal Sharing method. Furthermore, we also examine the relationship of these methods to the solution of the associated cooperative game.}
We conclude the methodological study by introducing a general family of segment allocation methods \textcolor{black}{that includes the three aforementioned rules. Finally, we evaluate the performance of these methods using a real-world dataset.}\end{abstract}


\end{frontmatter}

\textcolor{black}{
\paragraph{Highlights}
\begin{itemize}
 \item We revisit highway toll allocation problems from a game-theoretic perspective.
  \item We study two new allocation rules in this context using an axiomatic approach.
 \item We introduce a broad family of allocation rules.
 \item We apply our results to a real highway in Spain.
\end{itemize}}

\section{Introduction}

{Transportation is a key element of development policies in today's society. The role of infrastructure investment as an instrument of regional policy is reviewed in \cite{Asensio2023}. In \cite{Matas2018}, it is shown how, in the Spanish context, investment in infrastructure generates economic benefits when it alleviates pressure due to bottlenecks, connects strategic parts of the network, and reduces transport costs to the market. The problem of dealing with road maintenance and construction costs is a complex issue that can be addressed in various ways. Two fundamental questions emerge in this context: who should bear the financial responsibility and how much should they pay. \textcolor{black}{The primary response to the former} is either a toll system, where vehicles are charged for using a specific road, or budgetary allocation plans, where funds are collected from all taxpayers (regardless of their use of highways) through taxation, to finance road construction and maintenance. Additionally, alternative measures exist, such as the vignette system, which requires the payment of an annual fee for highway use. Regarding the second question,  namely how much should be charged for road use in the event that a toll system is adopted, a variety of different rates are established, with different charges being applied to different users depending on a range of factors. These include the user profile, the volume of traffic on the road, or the vehicle’s pollutant load.}

{As noted in \cite{Fiestras2011}, game theory has been successfully applied to the problem of cost allocation in the field of transportation. A notable example of this is the work of \cite{Littlechild1977}, where airport runway landing fee designs based on the Shapley value \citep{Shapley1953} and the nucleolus \citep{Schmeidler1969} are proposed. With the same purpose, in \cite{Vazquez1997} and \cite{Casas2003} two methods based on Owen's value \citep{Owen1977} and the $\tau$-value \citep{Tijs1981}, respectively, are introduced and axiomatically characterized, taking into account the aggregation of airplanes in airlines. More recently, \cite{Estan2021} and \cite{Babaei2022} have studied the allocation of a fixed cost between different cities participating in a line-based transport system, such as a tram or railway line. In the first case, three rules are proposed and characterized axiomatically, while in the second, the problem is analyzed through the lens of the well-known bankruptcy problem and the associated cooperative game. 
\cite{Ding2024} propose a two-stage hierarchical model to design a revenue-sharing rule that  ensures the stability of cooperation between different transportation service providers. The first stage solves the revenue sharing problem through a cooperative game, while the second stage identifies individualized resource allocation strategies through a non-cooperative game.} 

{The problem of distributing the cost of a highway among its users has also garnered interest in the literature. \cite{Dong2012} allocate these costs among all potential users and deal with the axiomatizations of a cost-sharing method (toll pricing method). \cite{Kuipers2013} use a similar model, but their focus lies on computing the Shapley value and the nucleolus of the so-called highway game that represents the problem. \cite{Gomez2024} assume that different classes of vehicles or travelers use the highway in such a way that their bargaining power can generate discounts. They consider the highway game with an a priori unions structure, study properties in this context of the Owen value, the $\tau$-value with a priori unions and a new solution that combines the $\tau$-value and the Shapley value, concluding with an illustrative real-world example.} 
 
Recently, \cite{Wu2024} introduced the problem of a highway toll allocation. In this model, the highway is divided into several segments, and a certain toll is levied for each trip along the highway, that is, for each group of consecutive segments. \cite{Wu2024} discuss how to distribute the total collected toll among these segments. The authors introduce three methods based on different toll-charging schemes. Specifically, when the toll depends on the number of highway segments used by a particular vehicle, they propose the Segments Equal Sharing (or SES) method, which allocates the toll obtained from any user equally over the segments used by that user. One reason why this method is appealing is that it coincides with the Shapley value \citep{Shapley1953} for the associated cooperative game.
\textcolor{black}{In addition, two axiomatic characterizations of the SES method are provided based on the fairness properties. In particular, the highway toll allocation problem can also be seen as an extension of the polluted river problem introduced in \cite{Ni2007}. \textcolor{black}{In this problem, a non-negative cost is associated with each segment of a polluted river, and the question is how to allocate the cost of cleaning the river among the agents located in each segment.} 
Note, however, that highway toll allocation problems deal with revenues rather than costs. On the other hand, the game induced by the toll policy can be related to two well-known classes of cooperative games with restricted cooperation.
One is \cite{Myerson1977}'s nonnegative communication line-graph game, in which the line-graph is defined by the consecutive highway segments. The other is the game with a permission structure \citep{Gilles1992}, in which the permission structure is given by the orientation of the highway in a single direction of travel.} 

In this paper, we revisit the problem of distributing the toll among the highway segments. {\textcolor{black}{In particular}, we \textcolor{black}{take on} the challenging task of deepening the search for properties of a method of distributing the tolls collected on a highway. First, following \cite{Thomson2015}, we consider the upper and lower limits of allocations as an essential part of addressing the fair allocation problem. In particular, as additional references, the definition of upper and lower bounds serves as the basis for the fairness criteria proposed in \cite{Fragnelli2010}, where, for the first time, the Baker-Thompson rule is characterized in the context of allocating the cost of a runway among the airplanes that use it. \textcolor{black}{Furthermore, it also appears in \cite{Tijs1986}, where allocations are introduced based on separable and non-separable costs. Secondly, a key property of the SES method is that it reflects the change in the allocation of two neighboring highway segments when the highway is blocked, causing all tolls collected from users traveling through those two segments to be lost.}
However, a case can be made for a property that reflects the effect of such a blockage on the two resulting components, just as \cite{Herings2008} in the context of cooperative games with constrained communication via graphs. These new approaches lead to two new methods for the highway toll allocation problem, providing two alternative ways to distribute the tolls: 
the Segments Proportional Sharing (SPS) method and the Segments Compensated Sharing (SCS) method. Both methods are axiomatically characterized and their relationship to two well-known solutions for the associated cooperative game is established, thus providing an early response to some of the open questions raised in \cite{Wu2024}. Despite the differences between the three methods considered throughout this paper, a common pattern emerges that allows us to define a family of assignment methods that generalizes them. This provides a range of potential distributions of the total collected toll, enabling the selection of the optimal method based on the specific criteria or available information. Finally, we illustrate the methods presented with a dataset taken from a real-world example. In this dataset, we incorporate a new comparison between the SES, SPS and SCS methods regarding the degree of equity in the distribution using classical metrics. Note that equity properties or solutions of an equitable nature are considered appropriate in different contexts, such as the allocation of the state in bankruptcy problems \citep{Thomson2015}, the sharing of costs of an irrigation ditch \citep{Aadland1998} or associated with the installation of an elevator \citep{Alonso2020}.} 

The remainder of this paper is organized as follows. Section \ref{sec:prel} introduces key concepts from cooperative game theory and formulates the highway toll allocation problem. Section \ref{sec:meth} defines the new methods and provides their corresponding interpretations. Sections \ref{sec:axSPS} and \ref{sec:axSCS} present the characterizations of the proposed methods, both from an axiomatic and a game-theoretical point of view. Section \ref{sec:family} introduces a general family of toll allocation methods and studies the core for the segment allocation game. Section \ref{sec:application} applies the methods considered in the paper to a real-world scenario. Some concluding remarks are deferred to Section \ref{sec:conclusions}. Finally, two appendices are included.

\section{Preliminaries}\label{sec:prel}

This section introduces fundamental concepts of cooperative game theory and establishes the notation used throughout the paper. In addition, we present the highway toll allocation problem as described in \cite{Wu2024}, along with the most relevant results related to our study.

\subsection{Cooperative game theory concepts}

A situation in which a finite set of players can generate certain payoffs by cooperation can be described by a cooperative game with
transferable utility, or simply, a TU-game. A TU-game is defined as a pair $(N, v)$ where $N$ is a set of players and $v: 2^{N} \longrightarrow \mathbb{R}$ is a characteristic function on $N$, satisfying $v(\emptyset)= 0$. For each subset $S$ of $N$ (a coalition), $v(S) \in \mathbb{R}$ represents the worth that the players in $S$ can obtain when they combine their resources, without the help of other players outside $S$. For the sake of simplicity,  a TU-game $(N, v)$ will be represented by $v$. We denote the collection of all TU-games in $N$ by $G^{N}$. 

A payoff vector for $v \in G^{N}$ is an $|N|$-dimensional vector $x \in \mathbb{R}^{N}$ assigning a payoff $x_{i} \in \mathbb{R}$ to each player $i \in N$. A solution for TU-games with player set $N$ is a function $f: G^{N} \longrightarrow \mathbb{R}^{N}$, which maps each TU-game to a payoff vector. If the solution always assigns a unique payoff, we refer to it as a single-valued solution or allocation rule. Otherwise, it is a set-valued solution concept. The most important set-valued solution concept is the core $C(N,v)$, proposed by \cite{Gillies1953}, defined as the set of efficient and stable allocations:
\vspace{-0.2cm}\begin{equation}
C(N,v) = \left\{x \in \mathbb{R}^{N}\,:\, \sum_{i\in N}x_i = v(N) \mbox{ and } \sum_{i\in S}x_i \geq v(S) \mbox{ for all } S \subset N\right\}.\vspace{-0.2 cm}
\label{eq:corecoopgame}
\end{equation}
Consequently, the core is the set of payoffs that distributes the worth of the grand coalition among the players and guarantees that no coalition receives less than what it can get on its own. 

One of the most famous allocation rules for TU-games is the Shapley value \citep{Shapley1953}, given by $Sh_{i}(N,v)= \sum_{S \subseteq N\backslash \left\{i\right\}}p(S)(v(S \cup \left\{i\right\}) - v(S))$ for all $i \in N$ and $v \in G^{N}$, where $p(S) =|S|!(|N| -| S| - 1)!/|N|!$. Compromise values are also well-known solution concepts that were thoroughly analyzed by \cite{TijsOtten1993}. For any game, a compromise value selects a unique allocation between two vectors: maximal and minimal payoffs for players. Different ways of measuring players' maximal and minimal payoffs yield various solution concepts. By considering upper and lower bounds of the core as the maximum and minimum aspiration for players, \cite{Tijs1981} proposed the $\tau$-value for transferable utility games, defined, if it exists, as the convex combination of the minimal rights vector, $m(N,v)$ and the utopia vector, $M(N,v)${\footnote{Recall that if $i\in N$, $M_i(N,v)=v(N)-v(N\backslash \{i\})$ and $m_i(N,v)=\max_{S\,:\,i\in S}\left(v(S)-\sum_{j\in S\backslash \{i\}}M_j(N,v)\right)$.}}: $\tau(N,v) = m(N,v) + \alpha \cdot \left(M(N,v) - m(N,v)\right)$, with $\alpha \in \mathbb{R}$ such that $\sum_{i \in N}\tau_{i}(N,v) = v(N)$. 

A TU-game implicitly assumes that each agent can directly communicate with every other agent. However, in real-world scenarios, players' communication may be restricted for several reasons. The need to model these communication constraints in the context of cooperative games led to the introduction of new theoretical models, commonly referred to as communication games. These games are represented by an undirected communication graph, as introduced by \cite{Myerson1977}.  An undirected graph is a pair $(N,g)$, where $N$ is a set of nodes and $g$ is a collection of edges, ie $g \subseteq \left\{\left\{i,j\right\}: i, j \in N, i \neq j\right\}$ is a collection of subsets of $N$ such that each element of $g$ contains precisely two distinct elements of $N$. The vertices in the graph normally represent the players in the game, and the edges symbolize the communication links between them. Node $i$ and node $j$ can communicate if $\left\{i, j\right\} \in g$. 

A TU-game combined with a communication structure characterized by a graph is called a graph game or communication game and is expressed by a triple $(N, v, g)$, consisting of a finite set of players ($N$), a characteristic function on the set of coalitions of players ($v$), and a set of edges (communication links) between the players ($g$). A solution for graph games is a mapping that assigns a set of payoff vectors to every graph game. The best-known single-valued solution for graph games is the Myerson value \citep{Myerson1977}, a solution that assigns to any graph game the Shapley value for the so-called graph-restricted game. However, \cite{Herings2008} introduced a new solution for the class of communication games, named the average tree solution (hereafter, it will be referred to as the AT solution). 

\cite{Suzuki2010} applies the AT solution on games with a line-graph structure. Formally, a graph $(N,g)$ with a finite set of nodes $N$ is called a line-graph if the set of links is characterized as $g = \left\{\left\{1, 2\right\},\left\{2, 3\right\},...,\left\{n-1, n\right\}\right\}$. 
In a line-graph, each node can communicate only with its neighbor(s), i.e, the node(s) located next to it. Following \cite{Suzuki2010}, we will denote by $L_{i}$ and $R_{i}$ the set of nodes located left to $i$ and right to $i$, respectively, for each node $i \in N$. In other words, $L_{i} = \left\{1,...,i-1\right\}$ and $R_{i} = \left\{i+1,...,n\right\}$. {Thus, for any line-graph game, $(N, v, g)$, \cite{Suzuki2010} considers the possible payoffs to player $i$ with respect to another player $j$,  based on the position of $i$ relative to $j$ in the line-graph. That is, $v(R_{i} \cup \left\{i\right\}) - v(R_{i})$, if $j \in L_{i}$; $ v(N) - v(L_{i}) - v(R_{i})$, if $i=j$; and $ v(L_{i} \cup \left\{i\right\}) - v(L_{i})$, if $j \in R_{i}$.} 

The AT solution, for any $i \in N$,  can then be expressed in terms of the previous amounts as follows. 
\begin{equation}
\text{AT}_{i}(v,g) = \frac{1}{n}\bigg((i-1)(v(R_{i} \cup \left\{i\right\}) - v(R_{i})) + (v(N) - v(L_{i}) - v(R_{i})) + (n-i)(v(L_{i} \cup \left\{i\right\}) - v(L_{i}))\bigg).
\label{eq:ATSuzuki}
\end{equation}
For more details on the AT solution, the reader is referred to \cite{Herings2008} and \cite{Suzuki2010}.

\subsection{Highway toll allocation problem}

In order to formulate the highway toll allocation problem, \cite{Wu2024} consider a one-way linear tolled highway which is divided into $n$ segments, {indexed according to the natural ordering from 1 to $n$}. The set of segments is denoted by $N = \left\{1, 2,...,n\right\}$, where $n=|N|$. Moreover, the entrance of segment $i \in N \backslash \left\{1\right\}$ is located at the same place as the exit of segment $i - 1$. For $h < k$, $[h,k]$ represents the path from the
entrance of $h$ to the exit of $k$, i.e, the sequence of consecutive segments $(h, h+1,...,k-1, k)$. The set of segments on the path $[h,k]$ is denoted by $\mathcal{P}([h, k])$, i.e., $\mathcal{P}([h, k])=\left\{h,...,k\right\}$. For technical convenience, an isolated segment is also treated as a path. The toll collected from all users entering at entrance $i$ and leaving at exit $j$ is symbolized by $t_{ij} \geq 0$. An $n\times n$-dimensional non-negative matrix $T$ is called a one-way toll matrix (or toll matrix for short) if $t_{ij} = 0$, for every $i >j$. The collection of all $n\times n$-dimensional toll matrices is denoted by $\mathcal{T}^{N}$. Finally, the set of all possible trips through the highway is given by $\left\{[h,k]: h,k \in N, h \leq k\right\}$. 

A highway toll allocation problem is a pair $(N, T)$ with $N =\left\{1, 2,...,n\right\}$ and $T \in \mathcal{T}^{N}$. A toll allocation for a problem $(N, T)$ is a non-negative vector $x \in \mathbb{R}^{N}_{+} $ assigning a share $x_{i} \in \mathbb{R}_{+}$ in the total toll to each segment $i \in N$. A toll allocation method (or method for short) is a mapping $f$ that assigns a toll allocation, $f(N, T)$, to each highway toll allocation problem $(N, T)$. Since we take the set of segments $N$ to be fixed, we often write a highway toll allocation problem, $(N, T)$, simply by its toll matrix, $T$. 

In the following, we present some results from \cite{Wu2024} that will be of interest along this paper. We firstly present the SES method for highway toll allocation problems, which distributes the toll revenue collected from each user equally over the segments used by that user.

\begin{definition}
\citep{Wu2024}. The Segments Equal Sharing method (shortly, the SES method) is given. for every $T\in \mathcal{T}^N$, by
\begin{equation*} 
f^{Se}_{i}(N,T) = \sum_{h=1}^{i}\sum_{k=i}^{n}\frac{t_{hk}}{k-h+1},\ {\rm for}\ {\rm every}\ i\in N.
\end{equation*}
\end{definition}

Once a toll allocation method has been formulated, it would be interesting to study the properties it verifies. We present below a set of natural axioms that any toll allocation method $f$ could satisfy, and which are, in particular, fulfilled by the SES method.

\textbf{Efficiency} A toll allocation method $f$ satisfies this property if, for every $T\in \mathcal{T}^N$, it holds that $\sum_{i\in N}f_i(N,T)=\sum_{i\in N}\sum_{j\in N}t_{ij}$.

\textbf{Inessential segment} A toll allocation method $f$ satisfies this property if, for every $T\in \mathcal{T}^N$ and for all $i\in N$ such that $t_{hk}=0$ for all $h\leq i\leq k$, then $f_i(N,T)=0$.

\textbf{Additivity} A toll allocation method $f$ satisfies this property if, for every pair $T, T'\in \mathcal{T}^N$, $f(N,T+T')=f(N,T)+f(N,T')$.

\textbf{Segment symmetry} A toll allocation method $f$ satisfies this property if, for every $T\in \mathcal{T}^N$ and every pair $i, j\in N$ such that $i,j\in \mathcal{P}([h,k])$ for all $h,k\in N$ with $t_{hk}> 0$, it holds that $f_i(N,T)=f_j(N,T)$.

\textbf{Toll fairness} A toll allocation method $f$ satisfies this property if, for every $T, T'\in \mathcal{T}^N$ and every $i\in N\backslash \{n\}$ such that $t'_{hk}=t_{hk}$, if $h\leq k\leq i$ or $i+1\leq h\leq k$ and $t'_{hk}=0$ if $h\leq i<k$ then $f_i(N,T)-f_i(N,T')=f_{i+1}(N,T)-f_{i+1}(N,T')$.

To present the following property, we consider the notion of sub-highway. A subset of consecutive segments $E\subseteq N$ is a sub-highway if $t_{ij}=0$ for every $i,j \in N$ with $\mathcal{P}([i,j])\cap E\neq \emptyset$ and $\{i,j\}\nsubseteq E$.

\textbf{Sub-highway efficiency} A toll allocation method $f$ satisfies this property if, for every $T \in \mathcal{T}^{N}$ and every sub-highway $E\subseteq N$, it holds that $\sum_{i\in E}f_i(N,T)=\sum_{h,k\in N:\, \mathcal{P}([h,k])\subseteq E}t_{hk}$.

Using these properties, \cite{Wu2024} obtain two axiomatic characterizations of the SES method.
\begin{theorem} \citep{Wu2024}.
\begin{enumerate}
\item The SES method is the only method that satisfies additivity, efficiency, the inessential segment property, and segment symmetry.
\item The SES method is the only method that satisfies toll fairness and sub-highway efficiency.
\end{enumerate}
\end{theorem}

A specific TU-game, the segments-allocation game, can be associated with the SES method. In this game, each coalition is assigned the toll collected from using only the segments included in that coalition.

\begin{definition} \citep{Wu2024}.
For every $T \in \mathcal{T}^{N}$, the segments allocation game $(N,v^{Se})$ is defined as
\begin{equation} 
v^{Se}(S) = \sum_{h, k \in N: \mathcal{P}([h,k]) \subseteq S}t_{hk},\ {\rm for}\ {\rm every}\ S\subseteq N.
\label{eq:segallgame}
\end{equation}
\end{definition}
 
\cite{Wu2024} also prove that every segment allocation game is a non-negative communication line-graph game 
and, moreover, every non-negative line-graph game is a segments allocation game. In addition, \cite{Wu2024} show that the allocation provided by the SES method coincides with the Shapley value for the associated segments allocation game.

\begin{theorem} \citep{Wu2024}.
Let $T \in \mathcal{T}^{N}$. Then $f^{Se}(N,T)= Sh(N, v^{Se})$.
\end{theorem}

\section{Toll allocation methods}\label{sec:meth}

This section considers two alternative toll allocation methods, hereinafter referred to as the Segments Proportional Sharing method and the Segments Compensated Sharing method. Both will be discussed in the following sections according to the outline of \cite{Wu2024} for the case of the SES method. Below, we also present the theoretical basis that justifies these methods, along with the notation required for their formal definition.

To motivate the usage of the first alternative, we introduce the following considerations. Fix $T\in \mathcal{T}^N$. The two-stage cost allocation method based on separable and non-separable costs (\citealp{Tijs1986}) can be adapted to $T$ leading to a distribution that takes into account the minimum and maximum revenue a highway segment in $T$ can achieve. 
Specifically, each segment $i\in N$ is first assigned its marginal or separable revenue $S_i(T)=t_{ii}$, which is a lower bound of segment $i$'s allocation. This quantity can be fairly ensured by yielding to the remaining segments the tolls collected on trips that make use of any segment in $N\backslash \{i\}$. For every highway toll allocation problem, it clearly satisfies that $\sum_{i\in N}t_{ii}\leq \sum_{i\in N}\sum_{j\in N}t_{ij}$. Secondly, the non-separable income, denoted by $NS(T)=\sum_{i\in N}\sum_{j\in N}t_{ij}-\sum_{i\in N}t_{ii}$, is distributed among all the segments. So, if $NS(T)>0$, a generic method could assign, to each $i$, the amount 
$S_i(T)+w_i \cdot NS(T)/(\sum_{j\in N}w_j)$, where $w_1,\dots,w_n$ are non-negative real numbers. In this setup, we propose to use the non-separable income involving segment $i$, denoted by $NS_i(T)=\sum_{h=1}^i\sum_{k=i}^nt_{hk}-t_{ii}$ for each $i\in N$. Note that $S_i(T)+NS_i(T)$ represents an upper bound on the allocation to segment $i$, as it is fair that the remaining segments can jointly secure the additional portion of the tolls earned since they come from trips not involving segment $i$. 

Thus, we introduce the Segments Proportional Sharing method that allocates, for every $T\in \mathcal{T}^N$, the tolls obtained from the users according to the previous ideas. 
{\begin{definition}
The Segments Proportional Sharing method (shortly, the SPS method) is given, for every $T\in \mathcal{T}^N$, by
\begin{equation} \label{eq:fSp}
f^{Sp}_{i}(N,T) = t_{ii} +{\displaystyle \frac{\sum_{h=1}^{i}\sum_{k=i}^{n}t_{hk}-t_{ii}}
{\sum_{j=1}^{n}\sum_{h=1}^{j}\sum_{k=j}^{n}t_{hk} - \sum_{j=1}^{n}t_{jj}}}
\cdot 
\left(\sum_{h=1}^{n}\sum_{k=h}^{n}t_{hk} - \sum_{j=1}^{n}t_{jj}\right),
\ {\rm for}\ {\rm every}\ i\in N.
\end{equation}
\end{definition}}

In the search for fairness and user acceptance of methods, alternative toll allocation procedures consider egalitarian methods with certain restrictions \citep{Aadland1998}. The following method, called Segments Compensated Sharing method, assigns to each segment a weighted average of the tolls of the different trips in which it is involved with the weight depending on the segment's position in the trip whether it is an entry or an exit.  Specifically, it compensates for the fact that a segment at the end of the highway potentially appears less frequently at the start of a trip and that a segment at the beginning that can be used less frequently as an exit.



Below, we present the Segments Compensated Sharing method as a toll allocation procedure for  every $T\in \mathcal{T}^N$.
\begin{definition}\label{def:AT1}
The Segments Compensated Sharing method (shortly, the SCS method) is given, for every $T\in \mathcal{T}^N$, by\vspace{-0.2cm}
\begin{equation}
f^{Sc}_{i}(N,T) = \frac{1}{n}\left((i-1)\sum_{k=i}^{n}t_{ik} + \sum_{h=1}^{i}\sum_{k=i}^{n}t_{hk} + (n-i)\sum_{h=1}^{i}t_{hi}\right), \ {\rm for}\ {\rm every}\ i\in N.
\label{eq:formorAT1}
\end{equation}
\end{definition}
In the previous expression,  $\sum_{k=i}^{n}t_{ik}$ is the toll collected from users entering at entrance $i$ of the highway, $\sum_{h=1}^{i}\sum_{k=i}^{n}t_{hk}$ is the total toll collected from all the highway trips that involve segment $i$ and $\sum_{h=1}^{i}t_{hi}$ is the toll collected from those trips where segment $i$ is used to exit the highway.
\begin{example}
Consider the toll allocation problem $\mathcal{T}^N$,  with $N = \{1, 2, 3\}$, and specified by  $t_{12}=t_{13} = 1$, and $t_{ij} = 0$, otherwise. Thus, we obtain  $f^{Se}(N,T) = (5/6, 5/6, 1/3)$,
$f^{Sp}(N,T)= (4/5, 4/5, 2/5)$, and $f^{Sc}(N,T)= (2/3, 1, 1/3)$.
\end{example}

\section{Characterization of the Segments Proportional Sharing method}\label{sec:axSPS}

In this section we theoretically study the SPS method. First, we characterize the SPS method from an axiomatic perspective in Subsection \ref{SPSproperties}. In Subsection \ref{SPS_games}, we prove, from a game-theoretic approach, that the allocation according to this method coincides with the $\tau$-value for the segments allocation game.

\subsection{Properties}\label{SPSproperties}

The first property is covariance, which states that if we change the unit of measurement for the tolls (for example, from euro to dollars) by multiplying the tolls by a given rate, the share allocated to each segment changes by the same amount. If, in addition, the toll of the trips that use only one segment increases by a constant amount, then the allocation of every segment will also increase in a proportionate manner. \par

\textbf{Covariance} A toll allocation method $f$ satisfies this property if, for every $T, T' \in \mathcal{T}^{N}$, $b \in \mathbb{R}_{+}$ and $a \in \mathbb{R}^{N}$, such that $t'_{ij} = b \cdot t_{ij}$ for $i \neq j$ and $t'_{ii} = b \cdot t_{ii} + a_{i}$, then $f_{i}(N,T') = bf_{i}(N,T) + a_{{i}}$, for all $i \in N$. \par

The second property we consider is the weighted segment symmetry property, which states that if all highway users use more than one segment on their trips, then the ratio between the assignments to each segment is precisely the ratio between the tolls collected from the trips using each of the segments. \par

\textbf{Weighted segment symmetry} A toll allocation method $f$ satisfies this property if, for every $T \in \mathcal{T}^{N}$ such $t_{ii}=0$ for all $i \in N$, then,  for all pair of essential segments\footnote{{A segment $i$ is said to be essential if $t_{hk}\neq 0$ for some $h,k\in N$ such that $h\leq i\leq k$.}} $i,j \in N$$$\frac{f_{i}(N,T)}{f_{j}(N,T)} = \frac{\sum_{h=1}^{i}\sum_{k=i}^{n}t_{hk}}{\sum_{h=1}^{j}\sum_{k=j}^{n}t_{hk}}.$$ 
Note that the property of weighted segment symmetry implies that segment symmetry is satisfied.

Below we provide a characterization of the SPS method based on efficiency, covariance, weighted segment symmetry and the inessential segment property.
\begin{theorem}\label{th:ch:SPS}
The SPS method is the only method that satisfies efficiency, the inessential segment property, the weighted segment symmetry property, and the covariance property.
\label{th:tauvalue}
\end{theorem}

An analysis of the logical independence of the axioms considered in Theorem \ref{th:ch:SPS} is in \ref{appA1}.

\subsection{Relation with the segments allocation game}\label{SPS_games}

The following result ensures that the SPS method, for any toll allocation problem $T \in \mathcal{T}^{N}$, coincides with the $\tau$-value for 
the associated segments allocation game $(N, v^{Se})$.

\begin{theorem}
Let $T \in \mathcal{T}^{N}$. Then $f^{Sp}(N,T) = \tau(N, v^{Se})$.
\end{theorem}

\section{Characterizations of the Segments Compensated Sharing method}\label{sec:axSCS}

In this section, we provide two characterizations of the SCS method, which are comparable to the two axiomatizations of the SES method presented in \cite{Wu2024}. Besides, we prove that the AT solution for the segments allocation game (which is a communication line-graph game) coincides with the allocation according to the SCS method.

\subsection{Properties}

The first property is weak segment symmetry, 
a weaker version of the segment symmetry property of \cite{Wu2024}.  
It states that in the absence of users that traverse the highway (except possibly some users who take the trip $[1,n]$), every segment should receive the same payment. Together with efficiency, it implies that each segment receives a revenue of $t_{1n}/n$, which means that the toll collected in the trip that uses all the highway segments is evenly distributed among them. 

\textbf{Weak  segment symmetry} A toll allocation method $f$ satisfies this property if, for every $T \in \mathcal{T}^{N}$, if $t_{hk}=0$ for all $h,k \in N$ such that $[h,k] \neq [1,n]$, then $f_{i}(N,T) = f_{j}(N,T)$, for all $i,j \in N$. \par

Before we formulate the next property, we need to present some notation. Note that every toll matrix $T \in \mathcal{T}^{N}$ can be decomposed as $T = \sum_{h=1}^{n}\sum_{k=h}^{n}t_{hk}\cdot {\delta^{hk}}$, where ${\delta^{hk}}$ is the matrix defined as $\delta^{hk}_{ij} = 1$, if $i=h, j=k$ and 0, otherwise, for each $h, k \in N$. The matrices ${\delta^{hk}}$ can also be thought of as toll matrices. Specifically, each $\delta^{hk}$ could represent the scenario in which only 1 unit is charged for the trip $[h,k]$ and only one user undertakes that trip. From now on, we will refer to these matrices as unitary toll matrices. 

We now introduce the indifference to individual extensions property. For any trip $[h,k]$, we consider two unitary toll matrices based on two possibilities of extending trips. Each of them includes an extra segment $j$, either placed just before $h$ or just after $k$. The indifference to individual extensions property means that the payoffs for the segments in $[h,k]$, except the ones directly connected to $j$, remain unchanged for both matrices. That is, adding a new segment to a trip only affects the segment it directly connects to.


\textbf{Indifference to individual extensions} A toll allocation $f$ satisfies this property if, for every $T\in \mathcal{T}^N$ and for every $h,k\in N$ with $h\leq k$, it holds that $f_{i}(N, \delta^{hk}) = f_{i}(N, \delta^{jk})$ for all $i \in \mathcal{P}([h,k])$ with $i\neq h$ and $j=h-1$, and that $f_{i}(N, \delta^{hk}) =f_{i}(N, \delta^{hj})$ for all $i \in \mathcal{P}([h,k])$ with $i\neq k$ and $j=k+1$.

Below, we present the linearity property, which encompasses both the additivity and the scale covariance properties. First, note that if $T\in \mathcal{T}^N$ and $b\in \mathbb{R}_+$, then $b\cdot T\in \mathcal{T}^N$ is a new toll matrix where $(b\cdot T)_{ij}=b\cdot t_{ij}$, for all $i,j\in N.$

\textbf{Linearity} A toll allocation method $f$ satisfies this property if, for every $T, T' \in \mathcal{T}^{N}$ and $b, b' \in \mathbb{R}_{+}$, $f(N,b\cdot T+b'\cdot T') = b\cdot f(N,T)+b'\cdot f(N,T')$. \par

Weak segment symmetry and indifference to individual extensions, along with linearity, efficiency and the inessential segment property, characterize the SCS method. Theorem \ref{th:ax1AT} proves such result.

\begin{theorem}\label{th:ax1AT}
The SCS method is the only method satisfying efficiency, linearity, the inessential segment property, weak segment symmetry and indifference to individual extensions. 
\end{theorem}

\ref{appA2} includes, also in this case, the analysis of the logical independence of the  axioms considered in Theorem \ref{th:ax1AT}.

Now, we provide a second characterization of the SCS method, which is comparable to the axiomatization of the SES method by means of sub-highway efficiency and toll fairness. We keep sub-highway efficiency as an axiom and consider an alternative fairness property, to be called toll component fairness. 

To justify the toll component fairness property, consider the highway toll allocation problem where the connection between two consecutive segments, $i$ and $i+1$, is ``blocked''. This means that no trip starting at or before entrance $i$ and ending at or after exit $i+1$ can be made. As a result, the tolls collected from all trips involving these two segments are lost, since blocking the connection disconnects not only segments $i$ and $i+1$, but also separates all segments before $i$ from those after $i+1$, creating two disconnected components. {Therefore, the losses associated with “blocking” a connection should be attributed to the two components, rather than to the two individual segments whose connection is blocked.}

\textbf{Toll component fairness} A toll allocation method $f$ satisfies this property if, for every $T, T' \in T^{N}$ and every $i \in N \backslash \left\{n\right\}$ such that $t'_{hk} = t_{hk}$ if $h \leq k \leq i$ or $i+1 \leq h \leq k$, and $t'_{hk}=0$ if $h \leq i < k$, it holds that $\frac{1}{i}\sum_{h \leq i}(f_{h}(N,T) - f_{h}(N,T')) = \frac{1}{n-i}\sum_{k>i}(f_{k}(N,T) - f_{k}(N,T'))$. \par

Toll component fairness states that blocking the connection between two segments yields for both resulting components the same average loss in payoff, where the average is taken over the segments in that component. 


Theorem \ref{th:atsol} presents the alternative axiomatization of the SCS method.
\begin{theorem}\label{th:atsol}
The SCS method is the only method satisfying toll component fairness and sub-highway efficiency.
\end{theorem}

\ref{appA2} provides an analysis of the logical independence of the  axioms considered in Theorem \ref{th:atsol}.

\subsection{Relation with the segments allocation game}

The following result ensures that the SCS method, for any toll allocation problem $T \in \mathcal{T}^{N}$, coincides with the AT solution for the associated segments allocation game. 

\begin{theorem}
Let $T \in \mathcal{T}^{N}$. Then $f^{Sc}(N,T) = AT(N, v^{Se})$.
\end{theorem}

\section{A general family of methods for highway toll allocation problems}\label{sec:family}

The aim of this section is to define a broader framework of toll allocation methods that encompasses those discussed in this paper. The shared features that the SPS, SCS, and SES methods present make this goal attainable. Among these shared features, we highlight the following two relevant issues:
\begin{itemize}
\item To determine the allocation of revenue to each segment, these three methods only consider those highway trips that involve such segment. 
\item Each considered method assigns a weight to the toll of each trip, $t_{hk}$, which can also be interpreted as a weight to the trip $[h,k]$. Such weight may be segment-dependent and it may vary across methods.
\end{itemize}

Using these ideas, we  introduce a general family of toll allocation methods for any $T\in \mathcal{T}^N$. First, we include some extra notation. We denote by $\mathcal{A}(N,T)$ the distribution of admissible weights for toll allocation in $T$ that is given by 
\begin{equation*}
	\mathcal{A}(N,T) = \left\{\alpha = \left(\alpha^{i}_{hk}(N,T)\right)^{i \in \mathcal{P}([h,k])}_{h,k\in N,\mbox{ }h\leq k}: \alpha^{i}_{hk}(N,T) \in \mathbb{R}_{+} \mbox{ for all }h,k\in N \mbox{ with }h\leq k \mbox{ and } i \in \mathcal{P}([h,k])  \right\}.
\end{equation*}
Thus, for any  toll allocation problem $T\in \mathcal{T}^N$ and for every $\alpha\in \mathcal{A}(N,T)$, a particular toll allocation for $T$ can be specified. Such method is denoted by $f^{\alpha}(N,T)$ and it is given, for every $i\in N$, by
\begin{equation}
f^{\alpha}_{i}(N,T) = \sum_{h=1}^{i}\sum_{k=i}^{n}\alpha^{i}_{hk}(N,T) \cdot t_{hk},
\label{eq:family}
\end{equation}for every $\alpha\in \mathcal{A}(N,T)$. Since the weights included in $\alpha$ are positive by definition, $f^{\alpha}(N,T)$ becomes a valid toll allocation method, i.e. $f^{\alpha}_{i}(N,T) \geq 0$ for every $i \in N$. However, the different nature of the weights, even under a non-linear dependence with $T$, does not ensure efficiency for the resulting toll allocation in a general framework. In practice, this fact implies that not all of the total toll collected along the highway is distributed among the segments.

Therefore, the family of toll allocation methods for every $T\in \mathcal{T}^N$ specified by the set of admissible weight distributions $\mathcal{A}(N,T)$ is given by the set $$\mathcal{F}(N,T) = \left\{f^{\alpha}(N,T) \in \mathbb{R}^{N}_{+}: \alpha \in \mathcal{A}(N,T)\right\}.$$ 

In the following, we firstly show the belonging of the three  methods considered along this paper (namely the SES, SPS and SCS methods) to the family $\mathcal{F}(N,T)$. Let $T$ be a highway toll allocation problem. Thus, specific distributions of weights $\alpha\in \mathcal{A}(N,T)$ give rise to different toll allocation methods. Fix $i \in N$, take $h,k\in N$ with $h\leq k$ such that $i \in \mathcal{P}([h,k])$ and consider the associated weight $\alpha^{i}_{hk}(N,T)$. Then, it holds that: 
\begin{itemize}
\item The SES method for $T$ is obtained when $\alpha^{i}_{hk}(N,T) = \displaystyle \frac{1}{k-h+1}$. 
\item The SPS method for $T$ results from considering 
\begin{equation*}
\alpha^{i}_{hk}(N,T) = \begin{cases} 1, \hspace{0.2cm}\text{if} \hspace{0.1cm} h=k; \vspace{-0.15 cm}\\
\beta, \hspace{0.2cm} \text{otherwise}\end{cases},\mbox{ where }\beta = \displaystyle \frac{\sum_{h=1}^{n}\sum_{k=h}^{n}t_{hk} - \sum_{j=1}^{n}t_{jj}}
{\sum_{j=1}^{n}\sum_{h=1}^{j}\sum_{k=j}^{n}t_{hk} - \sum_{j=1}^{n}t_{jj}}.
\end{equation*}

\item The SCS method for $T$ is obtained when 
\begin{equation*}
\alpha^{i}_{hk}(N,T) = \begin{cases} 1, \hspace{0.2cm} &\text{if} \hspace{0.1cm} h=k; \\ \frac{i}{n}, \hspace{0.2cm} &\text{if} \hspace{0.1cm} h = i \hspace{0.1cm} \text{and} \hspace{0.1cm} h \neq k; \\
\frac{n-i+1}{n}, \hspace{0.2cm} &\text{if} \hspace{0.1cm} k=i \hspace{0.1cm} \text{and} \hspace{0.1cm} h \neq k; \\
\frac{1}{n}, \hspace{0.2cm} &\text{otherwise}.
\end{cases}
\end{equation*}
\end{itemize}
In other words,  the Shapley value, the $\tau$-value and the AT solution for the segments allocations game belong to the same class of methods for highway toll allocation problems, despite their different nature. First, the weight distribution for the SES and SPS methods does not depend on the particular segment under consideration. However, this is not the case for the SCS method, leading to interesting interpretations. Specifically, if each segment is operated by a different company, it is reasonable to assume that a company might attribute greater importance to its segment for certain trips. For instance, due to higher traffic volumes, strategic positions (e.g., entrance or exit segments), or superior construction materials. Thus, the SCS method models an asymmetric scenario in which not all segments contribute equally or hold the same significance for a given trip. Moreover, although each segment within a trip is equally important under both the SES and SPS methods, these methods differ in how they weigh trips. The SES method assigns more relevance to shorter trips along the highway and less  to longer trips. In contrast, the SPS method values all trips equally, so each trip is equally important in calculating the share allocated to each segment. Note, however, that the ExES and the EnES methods, considered in \cite{Wu2024} for toll allocation problems, do not belong to the family $\mathcal{F}(N,T)$ either, as the allocation these methods assign to segment $i$ depends not only on trips involving this segment but also on trips where this segment is not included.

\subsection{On the stability of the general family of methods for the toll allocation problem} 

In this section we explore the relationship of the methods in $\mathcal{F}(N,T)$ for any highway toll allocation problem  and some well-known concepts in cooperative game theory, such as the core for a TU game (\ref{eq:corecoopgame}). 

The following proposition establishes that the core for the segments allocation game $(N, v^{Se})$ associated with any $T\in \mathcal{T}^N$ can be described by a particular subset of the family of methods $\mathcal{F}(N,T)$. These are specified by $T$-independent weights distributions $\alpha\in \mathcal{A}(N,T)$. A weight distribution $\alpha$ is said to be $T$-independent if $\alpha^{i}_{hk}(N,T)=\alpha^{i}_{hk}(N,T^\prime)$ for all  $T^{\prime}\in \mathcal{T}^N$ and for every $h,k\in N$ with $h\leq k$ and every $i \in \mathcal{P}([h,k])$. By technical convenience, under $T$-independence, we identify every $\alpha^{i}_{hk}(N,T)$, for every $h,k\in N$ with $h\leq k$ and for every $i \in \mathcal{P}([h,k])$, with $\alpha^i_{hk}$.  Furthermore, the distribution of weight has to guarantee that the toll collected on every trip is totally allocated among all the segments of such trip. Proposition \ref{prop:core} formalizes such result.

\begin{proposition}\label{prop:core} Let $T \in \mathcal{T}^{N}$. The core for the segments allocation game, $C(N,v^{Se})$, is given by
\begin{equation}
C(N,v^{Se}) = \left\{f^{\alpha}(N,T) \in \mathcal{F}(N,T): \alpha \mbox{ is  T-independent} \ \text{and } \ \sum_{i=h}^{k}\alpha^{i}_{hk} = 1 \ \text{for every} \ h, k \in N\right\}.
\label{eq:coresag}
\end{equation}
\end{proposition}

The next result states that the core is the set of those toll methods in $\mathcal{F}(N,T)$ satisfying efficiency, additivity and the inessential segment property. \par 

\begin{proposition}
Let $T\in \mathcal{T}^N$. A toll allocation $f(N,T) \in \mathbb{R}_{+}^{N}$ belongs to $C(N,v^{Se})$ if and only if $f(N,T)$ satisfies efficiency, additivity and the inessential segment property.
\label{prop:axcore}
\end{proposition}

By Propositions \ref{prop:core} and \ref{prop:axcore}, the SES method and the SCS method belong to the core for the segments allocation game. In contrast, the SPS method does not necessarily 
belong to the core, because the weights are toll-dependent {and, for every fixed trip, the weights of the segments involved do not add up to 1}. From a game-theoretic perspective, this means that for the segment allocation game, both the Shapley value and the AT solution always belong to the core, while the $\tau$-value does not. Example \ref{taunotinthecore} illustrates this statement.

\begin{example}\label{taunotinthecore}
	Take the toll allocation problem $T\in \mathcal{T}^N$,  with $N = \{1, 2, 3,4,5\}$, and specified by the upper triangular toll matrix $T$ such that $t_{11}=t_{15} = 1$, $t_{12}=5$, $t_{15}=1$, $t_{22}=1.50$, $t_{25}=0.02$, $t_{44}=2$, and $t_{ij} = 0.01$ for the remaining. Thus,  $f^{Sp}(N,T) = (3.401,3.917, 0.441, 2.427, 0.425)$ and it satisfies that
	 \[7.318 = 3.401+3.917=f_1^{Sp}(N,T)+f_2^{Sp}(N,T)<v^{Se}(\{1,2\})=7.5,\]
	and, hence, $f^{Sp}(N,T)\notin C(N,v^{Se})$.
\end{example}

Nevertheless, a necessary and sufficient condition for the inclusion of the SPS method in the core for $(N, v^{Se})$ can, in a certain sense, be formulated in accordance with the specific results presented in \cite{Driessen2013} for the $\tau$-value.

\begin{proposition}\label{tvinthecore}
  Let $T\in \mathcal{T}^N$.  The SPS method for $T$ belongs to $C(N,v^{Se})$ if and only if,  for every connected coalition of segments $S\subset N$,
  \begin{equation}\label{condbeta}
      \beta \geq  \frac{\underset{h, k \in N: \mathcal{P}([h,k]) \subseteq S}{\sum}t_{hk} - \underset{i\in S}{\sum}t_{ii}}
{\underset{i\in S}{\sum}\left(\underset{h, k \in N: i\in \mathcal{P}([h,k])}{\sum}t_{hk}-t_{ii}\right)}.
  \end{equation}
\end{proposition}

Clearly, the inequality in (\ref{condbeta}) does not hold in Example \ref{taunotinthecore}. However, it can easily be seen  that for highway toll allocation problems involving $|N|=3$ or $|N|=4$ segments, such condition is satisfied.

As \cite{Wu2024} stated, the tolls associated with segments $[h,k]$, with  $h,k\in N$ and $h\leq k$, correspond to the Harsanyi dividends of the segments allocation game associated to any $T\in \mathcal{T}^N$. Therefore, methods of the form in (\ref{eq:family}) can be interpreted as the possible ways in which the Harsanyi dividends of the trips involving a certain segment are distributed between the segments in those trips. Thus, the core for $(N,v^{Se})$ coincides with the Harsanyi set for such game.

\section{A practical application}\label{sec:application}

In this section, we apply the methods presented in this paper to a real-case scenario and compare them to the SES method. We will use the data from the Spanish highway AP68 considered in \cite{Wu2024} and in \cite{Kuipers2013} in a toll allocation framework. \par 

\noindent\textit{AP68 highway} \par

The highway AP68 connects Bilbao to Zaragoza, in the north of Spain, and has 23 entry/exit points along its route. The part of the highway between two adjacent points is naturally regarded as a segment. The general list of the segments of the AP68 highway, which form a line-graph, is shown in Table \ref{tab:segments}.

\begin{table}[h!]\begin{center}\resizebox{0.85\textwidth}{!}{
    \begin{tabular}{ll||ll||ll}
     \toprule
Index & Segment & Index & Segment& Index & Segment \\ \hline
1	&	 Bilbao - Arrigorriaga 	&	9	&	 Miranda  de  Ebro - Haro 	&	17	&	 Alfaro - E-AP15	 \\
2	&	 Arrigorriaga - Areta 	&	10	&	 Haro - Cenicero 	&	18	&	 E-AP15 - Tudela	 \\
3	&	 Areta - Llodio 	&	11	&	 Cenicero - Navarrete 	&	19	&	 Tudela - Gallur 	 \\
4	&	 Llodio - Ziorraga 	&	12	&	 Navarrete - Logroño 	&	20	&	 Gallur - A-272 	 \\
5	&	 Ziorraga - Altube 	&	13	&	 Logroño - Agoncillo 	&	21	&	 A-272 - A-275	 \\
6	&	 Altube - Subijana 	&	14	&	 Agoncillo - Lodosa 	&	22	&	 A-275 - Zaragoza	 \\
7	&	 Subijana - E-AP1 	&	15	&	 Lodosa - Calahorra	&		&		 \\
8	&	 E-AP1  - Miranda  de  Ebro 	&	16	&	 Calahorra - Alfaro	&		&		 \\     \bottomrule
    \end{tabular}}
\caption{Segments of the highway AP68.}
\label{tab:segments}\end{center}
\end{table}

Thus, we can formulate the associated highway toll allocation problem. 

\begin{table}[h!]\begin{center}\resizebox{0.9\textwidth}{!}{
    \begin{tabular}{llll||llll||llll}
     \toprule

$i$ & SES & SPS & SCS & $i$ & SES & SPS & SCS& $i$ & SES & SPS & SCS \\\hline
1	&	\makecell{30428.56 \\ (8.84 $\%$)} & \makecell{23261.38 \\ (6.76 $\%$)} &  \makecell{10331.42 \\ (3 $\%$)}	&	9 	&	\makecell{10006.06 \\ (2.91 $\%$)} & \makecell{13441.10 \\ (3.91 $\%$)} &	\makecell{15235.44 \\ (4.43 $\%$)} &	17 	& \makecell{8698.88 \\ (2.53 $\%$)}	 & \makecell{12027.74 \\ (3.49 $\%$)} & \makecell{7195.64 \\ (2.09 $\%$)}	 \\\hline
2	&	\makecell{29263.76 \\ (8.50 $\%$)} & \makecell{22096.58 \\ (6.42 $\%$)} &	\makecell{13923.71 \\ (4.05 $\%$)} &	10 	&	\makecell{9922.48 \\ (2.88 $\%$)} & \makecell{12960.35 \\ (3.77 $\%$)} & \makecell{11284.26 \\ (3.28 $\%$)} &	18	&	\makecell{15000.33 \\ (4.36 $\%$)} & \makecell{15556.73 \\ (4.52 $\%$)} &	 \makecell{26544.78 \\ (7.71 $\%$)} \\\hline
3	&	\makecell{26647.36 \\ (7.74 $\%$)} & \makecell{21523.22 \\ (6.25 $\%$)} & \makecell{15422.49 \\ (4.48 $\%$)} &	11	& \makecell{9682.34 \\ (2.81 $\%$)}	& \makecell{12740.88 \\ (3.7 $\%$)} & \makecell{8514.82 \\ (2.47 $\%$)} &	19	&	\makecell{15427.79 \\ (4.48 $\%$)} & \makecell{14863.03 \\ (4.32 $\%$)} &	\makecell{13537.09 \\ (3.93 $\%$)} \\\hline
4	& \makecell{25814.30 \\ (7.50 $\%$)} & \makecell{21256.30 \\ (6.18 $\%$)} & \makecell{9511.18 \\ (2.76 $\%$)}	&	12	&	\makecell{9551.17 \\ (2.78  $\%$)} & \makecell{12741.26 \\ (3.7  $\%$)} & \makecell{16758.75 \\ (4.87 $\%$)} &	20	&	\makecell{16962.08 \\ (4.93 $\%$)} & \makecell{15277.41  \\ (4.44 $\%$)} & \makecell{13672.06 \\ (3.97 $\%$)} \\\hline
5	&	\makecell{26390.88 \\ (7.67 $\%$)} & \makecell{21739.12 \\ (6.32 $\%$)} & 	\makecell{53338.43 \\ (15.5 $\%$)} &	13	&	\makecell{10413.31 \\ (3.03 $\%$)} & \makecell{13424.43  \\ (3.9$\%$)} & \makecell{16572.39 \\ (4.82 $\%$)} &	21	& \makecell{14755.64  \\ (4.29 $\%$)} & \makecell{13645.20 \\ (3.96 $\%$)} &	 \makecell{5733.30 \\ (1.67 $\%$)} \\\hline
6	& \makecell{14222.45 \\ (4.13 $\%$)} & \makecell{15368.54 \\ (4.47 $\%$)} & \makecell{7511.71 \\ (2.18 $\%$)}	&	14	& \makecell{10447.10 \\ (3.04 $\%$)}	& \makecell{12715.24 \\ (3.69 $\%$)} & \makecell{9686.86	\\ (2.81 $\%$)} &	22	&	\makecell{17025.39 \\ (4.95 $\%$)} & \makecell{15914.95 \\ (4.62 $\%$)} & \makecell{7930.37 \\ (2.3 $\%$)} \\\hline
7	&	\makecell{14177.34 \\ (4.12  $\%$)} & \makecell{15373.15 \\ (4.47 $\%$)} & \makecell{42878.96 \\ (12.46 $\%$)} &	15	&	\makecell{11020.85 \\ (3.2 $\%$)} & \makecell{13145.45 \\ (3.82 $\%$)} &	\makecell{10438.76 \\ (3.03 $\%$)} &		&	& &	 \\\cline{1-8}
8	& \makecell{9113.57 \\ (2.65 $\%$)}	& \makecell{12667.70 \\ (3.68 $\%$)} & \makecell{17697.8 \\ (5.14 $\%$)} &	16	& \makecell{9178.34 \\ (2.67 $\%$)}	& \makecell{12410.19 \\ (3.61 $\%$)} &	\makecell{10429.75 \\ (3.03 $\%$)} &		&	& &	 \\     \bottomrule
    \end{tabular}}
\caption{Toll allocation methods for the highway AP68.}
\label{tab:segmentsallocation}\end{center}
\end{table}

Table \ref{tab:segmentsallocation} shows the allocations specified per segment for the highway AP68 by the SES method, the SPS method and the SCS method, respectively, as well as the percentage of the total toll collected {that each allocation represents}. These results were obtained using the R software environment.\footnote{The \texttt{R} code built specifically for the determination of the toll allocations can be found at \url{https://github.com/PaulaSotoRodriguez/Highway_toll_allocation_methods}.} Figure \ref{fig:tollpoints} in \ref{App:dataset} completes this analysis through a graphical comparison. According to the allocation given for each of the methods mentioned in the paper, Table \ref{tab:ranking}  includes the top and bottom 3 segments.

\begin{table}[h!]\begin{center}\resizebox{0.9\textwidth}{!}{
 \begin{tabular}{l|ccc||l|ccc}
 \toprule
 Position & SES  & SPS  & SCS  & Position & SES  & SPS  & SCS\\\hline
 1 & Segment 1 & Segment 1 & Segment 5  & 20 & Segment 16 & Segment 8 & Segment 6  \\
 2 & Segment 2 & Segment 2 & Segment 7 & 21 & Segment 8 & Segment 16 & Segment 17 \\
 3 & Segment 3 & Segment 5 & Segment 18 & 22 & Segment 17 & Segment 17 & Segment 21\\
 \bottomrule
 \end{tabular}
 }
 \end{center}
 \caption{Ranking of the AP68 highway segments.}
\label{tab:ranking}
\end{table}

At first sight, the SES and SPS methods provides similar results, with slight variations in the positions of the same segments ($\rho_{Spearman}=0.947$). In particular, both methods allocate a larger toll to the first segment of the highway and a lower toll to segment 17. Segment 1 gets the most tolls mainly because the number of light vehicles entering from Segment 1 is much higher than other segments. Segment 17 gets the least tolls, not only because there are fewer light vehicles entering from the segment's entrance (ranked fifth from the bottom among all segments with entrances), but also because there are fewer light vehicles destined for Segment 17's exit (ranked last among all segments with exits). Thus, both the SES and SPS methods are consistent with the overall traffic flow situation on the AP68 highway. These methods are also almost identical in terms of the segments to which high tolls (above 15,000 euro) are allocated. Specifically, both methods allocate the highest tolls (above 25,000 euro in the case of the SES method and above 20,000 euro for the SPS method) to the first five segments of the highway. Moreover, for both these methods, the segments that receive the least amount are the ones located at the center of the highway (in other words, the most terminal segments receive the most). On the other hand, the SPS method distributes the toll equally among the segments. The SES method distributes the toll unevenly, forming groups/clusters of segments that receive more (such as the group of the first five segments and segments 18, 19, 20 and 22) and groups that receive much less. \par

For the SCS method, reader can check that the toll distribution differs from that of the SES and SPS methods. Note that, in this case, the ``central'' segments, such as segment 5, 7 or 18, receive the largest allocation, while the terminal segments, such as segment 21, receive the shortest proportion of the toll. This reflects an opposite behavior compared to the other two methods. In addition, the distribution of the SCS method is very uneven, with ``peaks'' of segments that receive a lot of money and others that receive very little. The range of toll amounts (difference between minimum and maximum) is very large, the maximum obtained being much higher than that obtained by the other two methods. It is particularly interesting that the SCS method allocates so much money to segment 5 (see Figure \ref{fig:tollpoints}). Although this segment is rarely used for entering the highway, it serves as the main exit for the majority of users. In fact, the most popular trip is the one that starts at entrance 1 and ends at exit 5, with 13061 users, so this trip collects the most tolls. It is notable that none of the other methods reflect this specific feature. Similarly, segment 7 is the second  most used segment for exiting the highway. In contrast, segment 18, located towards the end of the highway,  is mainly used to enter the highway. Indeed, it is the second most frequently used segment for this purpose, after segment 1, and therefore the third most favored by the SCS method. This method reflects nuances in the traffic flow on the highway that are not captured by the other methods. In particular, the SCS method reflects the individual behavior of those segments with special features, rather than considering groups of similar segments as is done in the SES method. \par

Table \ref{tab:corr} numerically summarizes these comparisons method by method. More in detail, we check the similarities between the SES and the SPS methods, as well as the differences with the allocation specified by the SCS method.
\begin{table}[h!]\begin{center}
    \begin{tabular}{l||ccc}
     \toprule
     & SES  &  SPS & SCS\\
     \hline
    SES	 &- &	0.989&	0.222\\
SPS	&0.947 & -&	0.284\\
SCS	&0.064& 	0.203&-\\
     \bottomrule
    \end{tabular}
\caption{Spearman’s correlation matrix (lower triangular matrix) and Pearson's correlation matrix (upper triangular matrix) for the toll allocations.}
\label{tab:corr}\end{center}
\end{table}

Finally, we calculate the corresponding Gini index \citep{Gini1912}, which is the most widely used measure to quantify the degree of socio-economic inequality, particularly in the context of benefit/income distribution. As Figure \ref{fig:correlations} also shows, both the SES and the SPS methods yield comparable results, indicating a positive and linear relationship between the two methods.  Given the similarity between these two methods, their relationship with the SCS method is also expected to be similar, as evidenced by the corresponding point clouds. We notice once more that the results of the SCS method are very different from the other two methods. \par

Figure \ref{fig:lorenzgini} shows the Lorenz curves and the Gini index for the three methods. The SPS method is the one that has the lowest Gini index, and thus represents the most equitable method. On the contrary, the SCS method has the highest Gini index. The SES method falls somewhere in between. Notably, although the SCS method is stable, it produces the least equitable outcome among the three methods. In contrast, the SPS method is the one that allocates the most equitable share of the collected toll.  \par

\section{Concluding remarks}\label{sec:conclusions}

Toll highways have become an essential method for alleviating congestion and supporting infrastructure development in today's society, highlighting the need for effective toll allocation strategies. In this paper we revisit the problem of toll sharing for a highway originally considered in \cite{Wu2024}. In such work,  the Segments Equal Sharing (SES) method was proposed as a toll-sharing mechanism, considering the distance traveled by vehicles as a basis. Using this framework, two new comparable methods are now innovatively considered: the Segments Proportional Sharing (SPS) method and the Segments Compensated Sharing (SCS) method. These proposals have been analyzed from a two-fold perspective. First, both toll methods are characterized by a set of natural properties tailored to a toll-allocation setting. Table \ref{tab:ax:cr} directly illustrates a comparative of the {properties} and the axiomatizations that the highway toll allocation methods discussed along this paper satisfy. Second, we have considered the relationship of both toll allocation procedures with well-known solution concepts from cooperative game theory, such as the $\tau$-value for TU-games and the AT solution for graph games. 

\begin{table}[h!]\begin{center}
\begin{tabular}{l||cccc}
\toprule
{Property}   & SES           &  SPS         & {SCS}       & Core \\\hline
Efficiency                 & $\checkmark^{2.1.1}$  & $\checkmark^{4.1}$  &$\checkmark^{5.1}$ & $\checkmark^{6.2}$ \\ 
Additivity                 & $\checkmark^{2.1.1}$  & -            & $\checkmark$ & $\checkmark^{6.2}$  \\
Linearity                  &  $\checkmark$ & -            &$\checkmark^{5.1}$ & $\checkmark$ \\
Covariance                 & $\checkmark$  & $\checkmark^{4.1}$ & $\checkmark$  & $\checkmark$ \\
Inessential segment property& $\checkmark^{2.1.1}$  & $\checkmark^{4.1}$ &$\checkmark^{5.1}$ & $\checkmark^{6.2}$ \\
Weak segment symmetry      & $\checkmark$  & $\checkmark$ & $\checkmark^{5.1}$ & - \\
Segment symmetry           & $\checkmark^{2.1.1}$  & $\checkmark$ & -  &  -    \\
Weighted segment symmetry  & -             & $\checkmark^{4.1}$ &- &  -       \\
Toll fairness              & $\checkmark^{2.1.2}$  & -            & -      &  -  \\
Toll component fairness    & -             & -            &$\checkmark^{5.2}$ & - \\
Sub-highway efficiency     & $\checkmark^{2.1.2}$  &   -          & $\checkmark^{5.2}$ & $\checkmark$ \\
Indifference to individual extensions  &   -    &     -   &  $\checkmark^{5.1}$ & - \\
     \bottomrule
    \end{tabular}
\caption{Overview of the properties and the axiomatizations for the highway toll allocation methods and the core for the segments allocation game.}
\label{tab:ax:cr}\end{center}
\end{table}

By their nature, we have extended these rules to a much larger family of toll methods, hereby introduced, that encompasses all those mentioned throughout this paper and (i) gives a a general way of constructing new methods for the highway toll allocation problem as well as (ii) new interpretations of the SES, SPS and SCS methods. We highlight that such a family does not necessarily provide efficient methods. In the context of the highway toll allocation problem, this fact has an interesting and realistic interpretation. For example, toll revenue is not always fully reinvested in highway maintenance or construction. Often, a portion is allocated to the highway operators or redirected to the state government. From a game-theoretic perspective, we have studied the belonging of the SES and the SCS methods to the core for the segments allocation game associated with any highway toll allocation problem. Finally, we apply the proposed toll allocation mechanisms on the Spanish highway AP68.
 
Below are some areas for further research and open issues that need to be addressed in the future, each heading in a different direction. First of all, one could consider extending the ideas and procedures of the toll allocation mechanisms presented here to the case of general graphs describing highways that are not necessarily linear. Similarly, extending segment allocation games to this broader context, as well as examining the relationship between the resulting allocations and specific solutions of these TU-games, would clearly be of interest. These aspects are essential for ensuring a fairer cost distribution, optimizing the management of transportation infrastructure, and promoting sustainable mobility in an ever-evolving society.

We also emphasize the close connection between the highway toll allocation problem and the museum pass problem introduced by \cite{Ginsburgh2001}. Specifically, any highway problem can be viewed as a museum pass problem with additional constraints: the museums must be visited in consecutive order, without skipping or backtracking. Furthermore, the Shapley value for the museum game \citep{Ginsburgh2003} shares the same underlying philosophy as the SES method, as it distributes the amount paid by a user equally among the museums they have visited. Exploring the interpretation of the SPS and SCS methods within the context of the museum problem and, conversely, examining how rules from the museum pass problem might be applied to the highway toll allocation problem could provide valuable insights. In this regard, we highlight the work of \cite{Martinez2023}, which generalizes the museum pass problem and characterizes several families of allocation rules. Investigating the relationship between these families and the family of toll allocation methods introduced in this paper presents another promising direction for future research.

\section*{Data and Code availability}
The \texttt{R} code of all the toll allocation methods presented in this paper, as well as the datasets used in the practical application, are available at \url{https://github.com/PaulaSotoRodriguez/Highway_toll_allocation_methods}.

\section*{Acknowledgments}

This work is part of the R+D+I project grant PID2021-124030NB-C32, funded by MCIN/AEI/10.13039/501
\noindent 100011033/, Spain and by “ERDF A way of making Europe”/EU. This research was also funded by Grupos de Referencia Competitiva ED431C 2021/24 from the Consellería de Cultura, Educación e Universidades, Xunta de Galicia, Spain.


\clearpage
\appendix

\section{Logical independence of the axioms}\label{App:indpendence}

In this appendix we study the logical independence of the axioms characterizing each of the toll allocation methods described along the paper.

\subsection{Logical independence of the axioms for the SPS method}\label{appA1}
\setcounter{theorem}{0} 
\renewcommand{\thetheorem}{4.1} 
\begin{theorem}
The SPS method is the only method satisfying efficiency, the inessential segment property, the weighted segment symmetry property and the covariance property.
\end{theorem}

Logical independence of the axioms in Theorem \ref{th:tauvalue} can be shown by the following alternative methods. 
\begin{enumerate}
\item Let $f$ be the method which allocates to each segment the toll collected from all the highway trips that use that segment.
\begin{equation*}
f_{i}(N,T) = \sum_{h=1}^{i}\sum_{k=i}^{n}t_{hk}, \hspace{0.2cm} \mbox{ for all } i \in N \hspace{0.2cm} \text{and for all } \hspace{0.2cm} T \in \mathcal{T}^{N}.
\end{equation*}

This method satisfies all properties with the exception of efficiency property.

\item Let $f$ be the combined method defined as 
\begin{equation*}
f(N,T) = \begin{cases} f'(N,T), \hspace{0.2cm} &\text{if} \hspace{0.2cm} T=T'; \\
f^{Sp}(N,T), \hspace{0.2cm} &\text{otherwise.}
\end{cases}
\end{equation*}
It suffices to take the toll allocation problem $(N,T')$, with $N=\{1, 2\}$, $t_{11}=1$, $t_{22}=2$, $t_{12}=0$ and   $f'(N,T')=(2,1)$, to state that this method verifies all properties except covariance.

\item Consider the SES method for toll allocation problems specified by
\begin{equation*}
f^{Se}_{i}(N,T) = \sum_{h=1}^{i}\sum_{k=i}^{n}\frac{t_{hk}}{k-h+1}, \hspace{0.2cm} \forall i \in N \hspace{0.2cm} \text{and} \hspace{0.2cm} T \in \mathcal{T}^{N}.
\end{equation*}
This method satisfies all properties except weighted segment symmetry.
\item  Let $f$ be the combined method defined as 
\begin{equation*}
f(N,T) = \begin{cases} f'(N,T), \hspace{0.2cm} &\text{if} \hspace{0.2cm} T\in \widetilde{T}^{N}; \\
f^{Sp}(N,T), \hspace{0.2cm} &\text{otherwise,}
\end{cases}
\end{equation*}

where $\widetilde{T}^{N}=\{T\in T^N\, :\, N=\{1,2,3\},\, t_{12}=t_{13}=0\ {\rm and}\ t_{23}>0\}$ and $f'_i(N,T)=t_{ii}+t_{23}/3$ for all $i\in N.$
This method satisfies all properties except the inessential segment property.
\end{enumerate}

\clearpage
\subsection{Logical independence of the axioms for the SCS method}\label{appA2}
\renewcommand{\thetheorem}{5.1} 
\begin{theorem}
The SCS method is the only method satisfying efficiency, linearity, the inessential segment property, weak segment symmetry and indifference to individual extensions. 
\end{theorem}

We will show that the axioms used in Theorem \ref{th:ax1AT} are independent.

\begin{enumerate}
    \item The SES method
    \begin{equation*}
    f_{i}(N,T) = \sum_{h=1}^{i}\sum_{k=i}^{n}\frac{t_{hk}}{k-h+1}, \hspace{0.2cm} \text{for every} \hspace{0.1cm} i \in N \hspace{0.2cm} \text{and for every} \hspace{0.2cm} T \in \mathcal{T}^{N}.
    \end{equation*}
    satisfies efficiency, the inessential segment property, weak segment symmetry and linearity, but fails indifference to individual extensions. This follows from the fact that for any unitary toll matrix, $\delta^{hk}$, the SES method assigns to every segment $i \in \mathcal{P}([h,k])$, $1/(k-h+1)$, while segments outside $\mathcal{P}([h,k])$ receive a payoff of zero. Thus, for trip $[j,k]$, with $|[j,k]| > |[h,k]|$ and $i \in \mathcal{P}([h,k])$ such that $i \neq h$, we will always have that $f_{i}(N, \delta^{hk}) \neq f_{i}(N, \delta^{jk})$.
    
    \item The method 
    \begin{equation*}
    f_{i}(N,T) = \frac{1}{n}\sum_{h=1}^{n}\sum_{k=h}^{n}t_{hk}, \hspace{0.2cm} \text{for every} \hspace{0.1cm} i \in N \hspace{0.2cm} \text{and  for every} \hspace{0.2cm} T \in \mathcal{T}^{N},
    \end{equation*}
    satisfies efficiency, weak segment symmetry, linearity and indifference to individual extensions, but does not verify the inessential segment property.
    
    \item The method
    \begin{equation*}
    f_{i}(N,T) = 0, \hspace{0.2cm} \text{for every} \hspace{0.1cm} i \in N \hspace{0.2cm} \text{and  for every} \hspace{0.2cm} T \in \mathcal{T}^{N},
    \end{equation*}
    satisfies weak segment symmetry, linearity, the inessential segment property and indifference to individual extensions, but fails efficiency.
    
    \item The method 
    \begin{equation*}
    f_{i}(N,T) = \sum_{k=i}^{n}t_{ik}, \hspace{0.2cm} \text{for every} \hspace{0.1cm} i \in N \hspace{0.2cm} \text{and  for every} \hspace{0.2cm} T \in \mathcal{T}^{N}, 
    \end{equation*}
    satisfies efficiency, the inessential segment property, linearity and indifference to individual extensions, but fails weak segment symmetry.

    \item The method \begin{equation*} f_{i}(N,T) = 
    \begin{cases}
    f_{i}^{Sc}(N,T), \hspace{0.2cm} \text{for} \hspace{0.2cm} T = \delta^{hk}, \hspace{0.1cm} h, k \in N;  \\
    f_{i}^{Se}(N,T), \hspace{0.2cm} \text{otherwise}, 
    \end{cases}
    \end{equation*}
    for every $i \in N$ and  for every $T \in \mathcal{T}^{N}$, satisfies efficiency, the inessential segment property, indifference to individual extensions and weak segment symmetry, yet fails linearity (take $T = \delta^{hk}$ and $T' = \delta^{h'k'}$, then $T + T'$ is not a unitary toll matrix, so $f(N, T+T') = f^{Se}(N, T+T')$, whereas $f(N,T) = f^{Sc}(N,T)$ and $f(N,T') = f^{Sc}(N,T')$).
\end{enumerate}

\renewcommand{\thetheorem}{5.2} 
\begin{theorem}
The SCS method is the only method satisfying toll component fairness and sub-highway efficiency. 
\end{theorem}

We show logical independence of the properties stated in Theorem \ref{th:atsol} by presenting two alternatives methods for highway toll allocation problems $T\in \mathcal{T}^N$.
\begin{enumerate}
\item Let $f$ be the method  in which no toll is allocated to any segment,
\begin{equation*}
f_{i}(N,T) = 0, \hspace{0.2cm} \text{for every} \hspace{0.1cm} i \in N \hspace{0.2cm} \text{and} \hspace{0.2cm} T \in \mathcal{T}^{N}.
\end{equation*}

Thus, $f$ satisfies toll component fairness, but not sub-highway efficiency.

\item Let $f$ be the SES method  
\begin{equation*}
f_{i}(N,T) = \sum_{h=1}^{i}\sum_{k=i}^{n}\frac{t_{hk}}{k-h+1}, \hspace{0.2cm} \text{for every} \hspace{0.1cm} i \in N \hspace{0.2cm} \text{and} \hspace{0.2cm} T \in \mathcal{T}^{N}.
\end{equation*}

This method satisfies sub-highway efficiency, but not toll component fairness.
\end{enumerate}

\clearpage
\section{The AP68 highway}\label{App:dataset}

In this appendix we include some graphical results that complete the ones shown in Section \ref{sec:application}.



\begin{figure}[H]
    \centering
    \scalebox{0.6}{
    \includegraphics{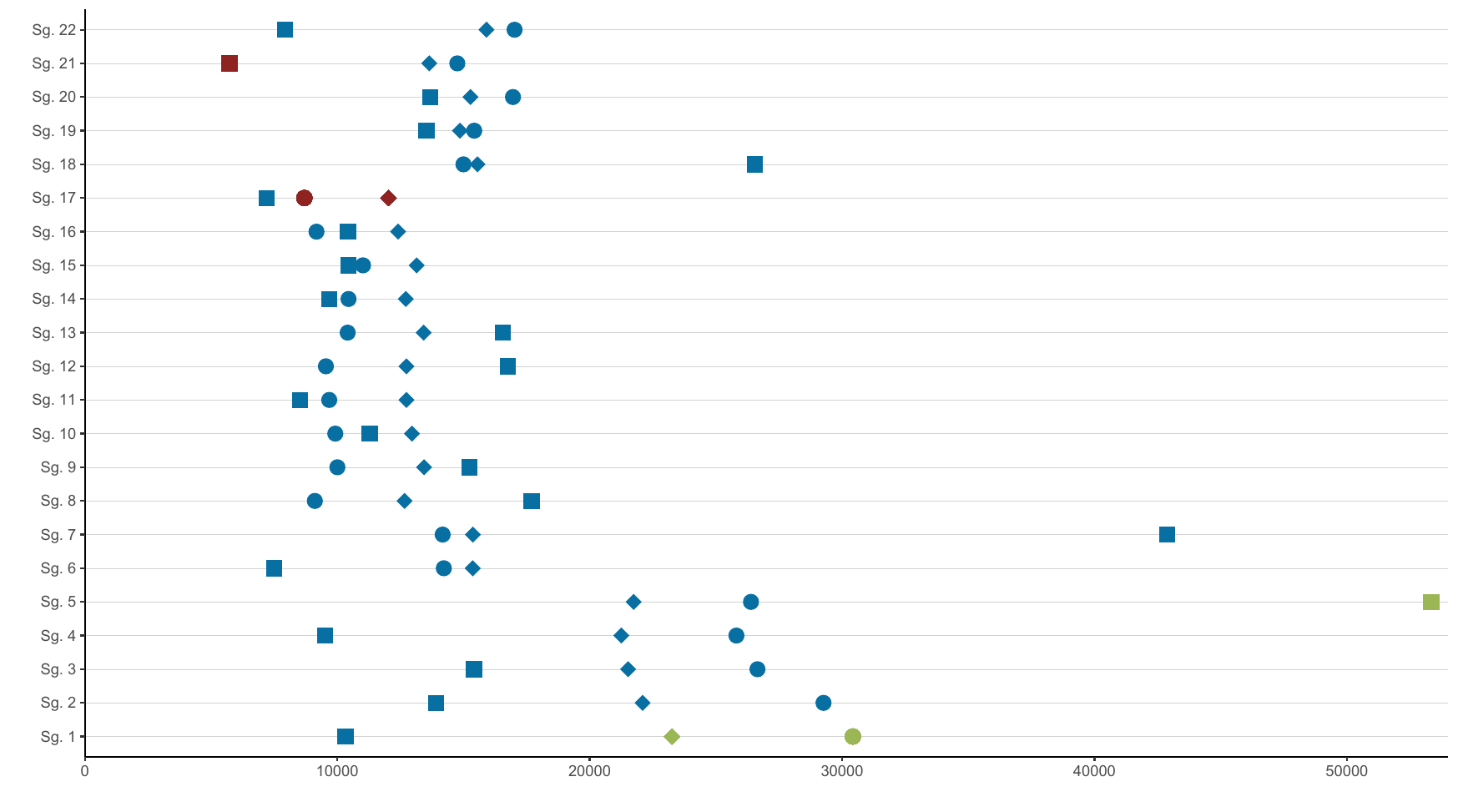}}
    \vspace{-0.25 cm}\caption{Toll allocation provided by the SES method (circle), SPS method (diamond) and SCS method (square) for the AP68 highway segments. The segment with the highest toll is marked in green, while the segment with the lowest toll is shown in red.}
    \label{fig:tollpoints}
\end{figure}

Once the  the SES method, the SPS method and the SCS method were applied for the case of the  highway AP68, the results were included in Table \ref{tab:segmentsallocation}. Figure \ref{fig:tollpoints} completes this analysis from a graphical perspective.

\begin{figure}[H]
    \centering
    \scalebox{0.6}{
    \includegraphics{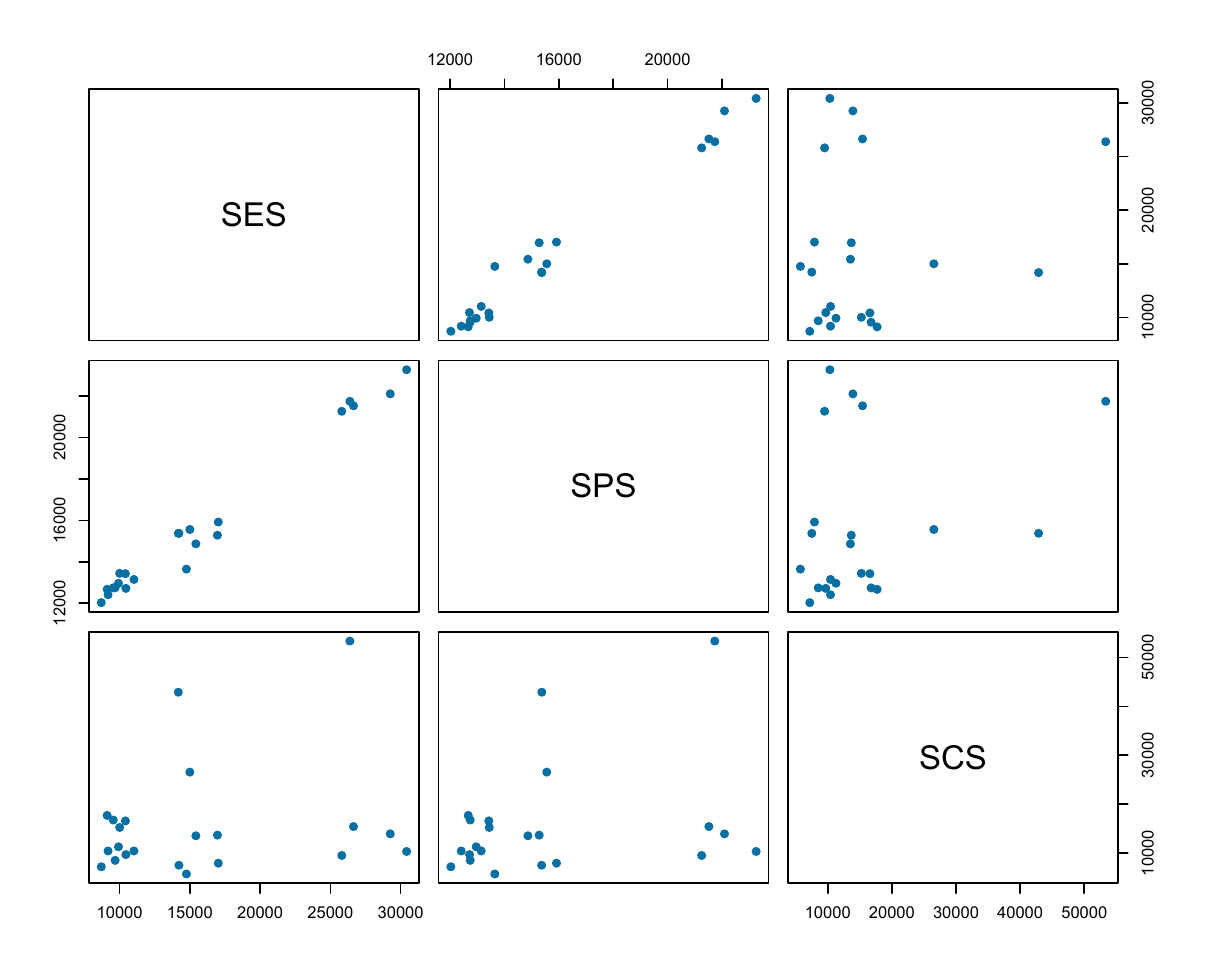}}
    \caption{Correlations between the SES, SPS and SCS methods.}
    \label{fig:correlations}
\end{figure}

In parallel to the numerical analysis of the correlations, Figure \ref{fig:correlations} shows the corresponding scatter plots between pairs of toll allocation methods.

\begin{figure}[H]
    \centering
    \scalebox{0.6}{
    \includegraphics{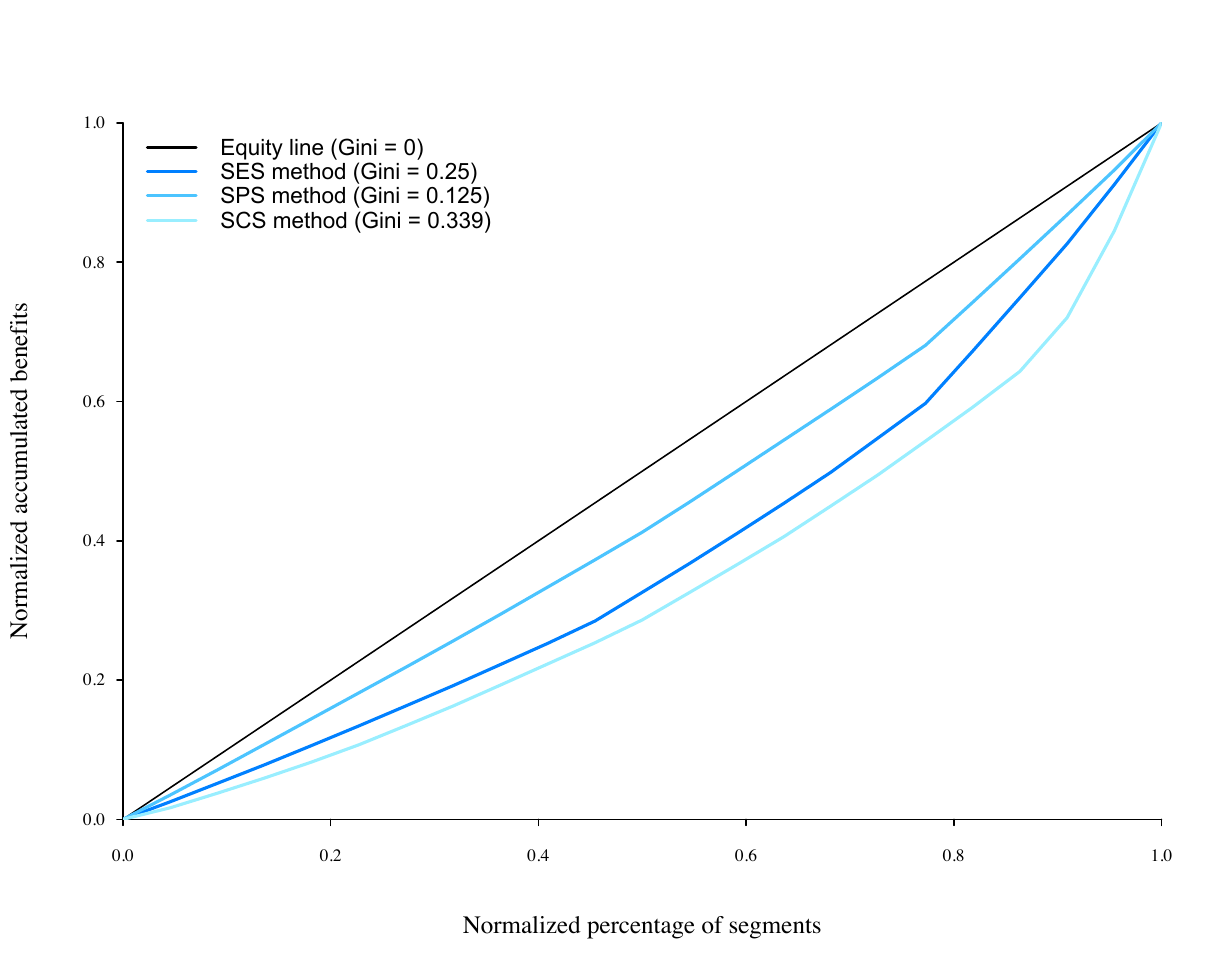}}
    \caption{Lorenz curves and Gini index for the SES, SPS and SCS methods.}
    \label{fig:lorenzgini}
\end{figure}

Finally, Figure \ref{fig:lorenzgini} represents the associated Lorenz curves for each of the toll methods considered.

\end{document}